\documentclass{amsart}

\usepackage[dvips]{color}
\usepackage{amsfonts}
\usepackage{multicol}
\usepackage{graphicx}
\usepackage{amsmath}
\usepackage{amssymb}
\usepackage{amsthm}
\usepackage{mathrsfs}
\newtheorem{theorem}{Theorem}[section]

\newtheorem{proposition}[theorem]{Proposition}

\theoremstyle{definition}

\theoremstyle{remark}

\numberwithin{equation}{section}

\newfam\msbfam
\font\tenmsb=msbm10  \textfont\msbfam=\tenmsb
\font\sevenmsb=msbm7  \scriptfont\msbfam=\sevenmsb
\font\fivemsb=msbm5    \scriptscriptfont\msbfam=\fivemsb
\def\Bbb{\fam\msbfam \tenmsb}

\def\CC{{\Bbb C}}

\newfam\msbbfam
\font\tenmsbb=msbm10  scaled \magstep1 \textfont\msbbfam=\tenmsbb
\font\sevenmsbb=msbm7  scaled \magstep1 \scriptfont\msbbfam=\sevenmsbb
\font\fivemsbb=msbm5    scaled \magstep1 \scriptscriptfont\msbbfam=\fivemsbb

\begin{document}

\title{Inner functions on a general type of domains}

\author{Baili Min}
\address{Department of Mathematics, Washington University in St.\,Louis, Saint Louis, MO 63130}
\email{minbaili@math.wustl.edu}

\date{\today}

\thanks{I would like to thank my advisor Prof.~Steven Krantz, who encouraged and directed to write this paper.}
\begin{abstract}
It is known that inner functions exist on strongly pseudoconvex domains. In this paper we will show that they exist on a more general type of domains, including some domains of finite type.
\end{abstract}
\maketitle
\section{Background and Introduction}
If $\Omega$ is a domain in $\CC^n$ and $f: \Omega \to \CC$ is a bounded holomorphic function such that for almost every $\zeta \in \partial \Omega$, the radial boundary limit $f^*(\zeta)$ exists and $|f^*(\zeta)|=1$, then we say that $f$ is an inner function. 
In this paper, we do not consider the trivial example of constant functions.

This subject has been fully studied in the case of a single variable. Illustrative examples on the unit disc are the Blaschke products with prescribed zeros $\{\alpha_i\}$ that satisfy the Blaschke condition:
\begin{displaymath}
B(z)=e^{i\theta}z^k\prod\frac{|\alpha_i|}{\alpha_i}\cdot\frac{\alpha_i-z}{1-\overline{\alpha}_i z}, \ \  |c|=1.
\end{displaymath}
Other examples of zero-free functions are
\begin{displaymath}
G(z)=\exp \Big\{-\int_{|\zeta|=1}\frac{\zeta+z}{\zeta-z}\,d\mu(\zeta)\Big\}
\end{displaymath}
where $\mu$ is a positive Borel measure on the unit circle. In fact, every inner function is a product of those two.

When it comes to the case of several variables, the problem is more complicated. At one time, the existence of such inner functions was doubted , and it was found that, even if an inner function existed, it had some undesirable properties, such as being discontinuous at every boundary point of the unit ball in $\CC^2$. For more discussion, please see \cite{R1} and \cite{K}. Later, a turnaround considerable the attention: inner functions were constructed for the unit ball, in $\CC^n$. For more on this, we refer the reader to the work of A.~Aleksandrov \cite{A1}, M.~Hakim and N.~Sibony \cite{HS}, and E.~L{\o}w \cite{LOW1}. Additionally, in \cite{LOW2}, L{\o}w showed that inner functions exist for strongly pseudoconvex domains. Their work uses various methods and tools, including Ryll-Wojtaszczyk polynomials, a method developed by Aleksandrov in \cite{A2} as an alternative approach to construct inner functions for the unit ball. Later, W.~Rudin wrote a book, \cite{R2}, on this method, and provided many other applications.
 
We keep asking ourselves, can we explore more general domains?

In this paper, inspired by Rudin's summarization in \cite{R2}, we present some results on domains that are similar to the unit ball insomuch as both are complex manifolds and there is a ramified holomorphic map between them. Although we present our results in the context of $\CC^2$, our method is generalizable to higher dimensions. Our principal results are Theorem \eqref{existencethm} and Proposition \eqref{existencepro}, which establish the existence of inner function on these domains.

\section{Integral Formulas}
Let $B$ be the unit ball in $\mathbb{C}^2$ and $S$ be its boundary.

Suppose: $f: \mathbb{C}^2 \to \mathbb{C}^2$ holomorphic, $M \subset \mathbb{C}^2$ is a compact manifold with smooth boundary such that $f(M)=B$, $f(\partial M)= S$.

Also suppose that $f=(f_1, f_2): \overline{M} \to \overline{B}$ is a finite ramified covering, with $N$ sheets. 
Let $Z=\{(z,w)\in \overline{M}: \text{either $z$ is a zero of $f_1$ or $w$ is a zero of $f_2$}\}$.
Denote $Z_1= Z \cap \partial M$, $Z_2=f(Z) \cap S$. 

We can then think of $\{U_1, \ldots, U_n\}$ is an open cover of $S \setminus Z_2$ such that for each $U_i$, $f^{-1}(U_i)$ contains $N$ disjoint components $V_{i1}, \ldots, V_{iN}$ such that each component is biholomorphic to $U_i$, and $\{V_{i1}, \ldots, V_{iN} \}_{i}$ covers $\partial M \setminus Z_2$.

Then we choose $\{W_1, \ldots, W_n\}$ such that $W_i \subset V_i$, $\cup_i^n W_i = S - Z_2$, and for each pair $(W_i, W_j)$, the interiors of $W_i$ and $W_j$ intersect at an empty set level.  If we let $\{X_{i1}, \ldots, X_{iN} \}=f^{-1}(W_i)$, then we know that $\cup_j^N\cup_i^n X_{ij} = \partial M - Z_1$ and for each $k$, the intersection of the interiors of $X_{ik}$ and $X_{jk}$ is an empty set.

This case is of interest because it can furnish domains that are more general than strongly pseudocovex ones. For example, for any positive integer $p$, the map $f_p(z_1, z_2)=(z_1, z_2^{2p})$ relates the domain $M_p=\{(z_1, z_2)\in \mathbb{C}^2: |z_1|^2+|z_2|^{2p}<1\}$, which is not strongly pseudoconvex but still of finite type, to the unit ball. There are many references for the notion of finite type, including \cite{C} by D.~Catlin.

For $z, w \in \partial M$ sufficiently close to one another, define 
\begin{equation}
d_M(z,w)=\sqrt{1-\big|\langle f(z), f(w) \rangle\big|^2},
\end{equation}
we can then introduce the open boundary ball, centered at $\omega \in \partial M$, with radius $0<r<1$:
\begin{equation}
E_M(\omega, r)=\{\zeta \in \partial M: d_M(\omega, \zeta)<r\}.
\end{equation}

We wish to introduce a measure over $\partial M$. The way to do this is to relate it to a measure on $S$, the boundary of the unit ball in $\mathbb{C}^2$. 

Specifically, for $\eta \in S$, put
\begin{equation}
E(\eta,r)=\{\xi \in S: d(\eta, \xi)=\sqrt{1-\big|\langle \eta, \xi \rangle\big|^2}<r\}, 
\end{equation}
and let $\sigma$ be the unique rotation-invariant probability measure on $S$.

With the map $f: \partial M \to S$, we immediately have the relations
\begin{equation}
d_M(z,w)=d(f(z), f(w)),
\end{equation}
and 
\begin{equation}
\label{boundaryballrelation}
f(E_M(\omega,r))=E(f(\omega),r).
\end{equation}
This inspires us to define a measure $\sigma_M$ by
\begin{equation}
\sigma_M(E_M(\omega, r))=\sigma(f(E_M(\omega,r)))=\sigma(E(f(\omega),r)),
\end{equation}
for $r$ sufficiently small.

Although here, $\sigma_M$ is defined for balls, it also works for open subsets of $\partial M$. We have the relation
\begin{equation}
\sigma_M=\sigma \circ f.
\end{equation}

The study of $E(\eta,\delta)$ shows that $\sigma(E(\eta,\delta))=\delta^2$. Consequently, we have
\begin{equation}
\sigma_M(E_M(\omega, r))=r^2.
\label{boundaryballarea}
\end{equation}

Now we turn our attention to the relation between the integrals over $S$ and those over $\partial M$.

For $j=1, 2, \ldots, N$, 
\begin{equation}
\bigcup_i^n \int_{X_{ij}} F\big(f(w)\big) \,d\sigma_M(w)=\bigcup_i^n \int_{W_i}F(z)\,d\sigma(z)=\int_{S-Z_2}F(z)\,d\sigma(z).
\end{equation}

Since $\sigma(Z_2)=0$, we have the equation
\begin{equation}
\bigcup_i^n \int_{X_{ij}} F\big(f(w)\big) \,d\sigma_M(w)=\int_{S}F(z)\,d\sigma(z).
\end{equation}

Therefore, since $\sigma_M(Z_1)=0$, we know the relation between these two kinds of integrals is:
\begin{align}
\int_{M}F\big(f(w)\big) \,d\sigma_M(w)&=\int_{M-Z_1}F\big(f(w)\big) \,d\sigma_M(w) \\ \nonumber
&=N\cdot\Big(\bigcup_i^n \int_{X_{ij}} F\big(f(w)\big) \,d\sigma_M(w)\Big)  \\  \nonumber
&=N\int_{S}F(z)\,d\sigma(z).
\end{align}

Now, let $\alpha=(\alpha_1, \alpha_2)$ be a multi-index. 
We use the notation $z^{\alpha}=z^{(\alpha_1, \alpha_2)}=z_1^{\alpha_1}z_2^{\alpha_2}$.

Since, as shown in Rudin's book, for $\alpha \neq \beta$,
\begin{equation}
\int_{S}z^{\alpha}\overline{z}^{\beta}\,d\sigma(z)=0,
\end{equation}
and
\begin{equation}
\int_{S}z^{\alpha}\overline{z}^{\alpha}\,d\sigma(z)=\frac{\alpha!}{(1+|\alpha|)!},
\end{equation}
we can then obtain similar equations for the case $\partial M$:
\begin{equation}
\label{zerointegral}
\int_{\partial M}f(z)^{\alpha}\overline{f(z)}^{\beta}\,d\sigma_M(z)=0, \alpha \neq \beta,
\end{equation}
and
\begin{equation}
\label{nonzerointegral}
\int_{\partial M}f(z)^{\alpha}\overline{f(z)}^{\alpha}\,d\sigma_M(z)=\frac{N\alpha!}{(1+|\alpha|)!}.
\end{equation}

\section{Boundary Balls}

Without too much difficulty, we can prove that, for $z, w$ and $u \in \partial M_2$, we have a triangle inequality, that is, there exists a positive constant $C_1$ such that
\begin{equation}
d_M(z, w) \leqslant C_1 (d_M(z, u)+d_M(u,w)).
\end{equation}

With the help of this general result, or, more simply, from the definition of the area of boundary balls as shown in \eqref{boundaryballarea}, along with what we established  in \eqref{boundaryballrelation}, we can check that:

1. $0<\sigma_M(E_M(\omega, r))<\infty$;

2. There exists $C_2>0$ such that $\sigma_M(E_M(\omega, 2r))\leqslant C_2 \sigma_M(E_M(\omega, r))$;

3. There exists $C_3>0$ such that if $E_M(\omega, r)\cap E_M(\zeta, s)=\emptyset$ and $r \leqslant s$, then $E_M(\omega, r) \subseteq E_M(\zeta, C_3 s)$.

This gives us a geometric result:
\begin{theorem}
For $r>0$, there exists a maximal set $\{\omega_1, \ldots, \omega_K\} \subset \partial M$ with respect to having the balls $E_M(\omega_j, r)$ pairwise disjoint but $\cup_{j=1}^K E_M(\omega_j, 2r)$ cover $\partial M$.
\end{theorem}

As a consequence of this, we have
\begin{equation}
C_4=\sigma_M(\partial M) \leqslant K \sigma_M(E_M(\omega_j, 2r))=K(2r)^2=4Kr^2,
\end{equation}
which implies that
\begin{equation}
K \geqslant \frac{C_4}{4r^2}.
\end{equation}

\section{Construction of $f$-Polynomials}
Let $r_1, \ldots, r_K$ be Rademacher functions, and define
\begin{equation}
Q_t \circ f(z)=\sum_{j=1}^K r_j(t)\langle f(z),f(\omega_j) \rangle^k,
\end{equation}
where $k \in \mathbb{N}$.

We first notice that $Q_t( \lambda f_1,  \lambda f_2)=\lambda^k Q_t \circ f$. Now, take $r=1/\sqrt{k}$. We wish to find bounds for $|Q_t \circ f|$.

\subsection{Lower Bounds}

We first calculate that
\begin{align}
\big|\langle f(\zeta), f(\omega_j) \rangle\big|^{2k}&=\big(f_1(\zeta)\overline{f_1(\omega_j)}+f_2(\zeta)\overline{f_2(\omega_j)}\big)^k 
\big(\overline{f_1(\zeta)}f_1(\omega_j)+\overline{f_2(\zeta)}f_2(\omega_j)\big)^k\\ \nonumber
                                                           &=\sum_{0 \leqslant i,s \leqslant k} D_{i,s}f(\zeta)^{(i,k-i)}\overline{f(\zeta)}^{(s,k-s)},
\end{align}
where
\begin{equation}
D_{i,s}={{k}\choose{i}}^2f(\omega)^{(s,k-s)}\overline{f(\omega)}^{(i,k-i)}.
\end{equation}

We note that $|f_1(\omega)|^2+|f_2(\omega)|^2=1$.
Using \eqref{zerointegral} and \eqref{nonzerointegral}, we have the equations
\begin{align}
\int_{\partial M}\big|\langle f(\zeta), f(\omega_j) \rangle\big|^{2k}\,d\sigma_M(\zeta)&= \int_{\partial M}\sum_{0 \leqslant i,s \leqslant k} D_{i,s}f(\zeta)^{(i,k-i)}\overline{f(\zeta)}^{(s,k-s)}\,d\sigma_M(\zeta) \\ \nonumber
                                                              &=\sum_{i=0}^k D_{i,i} \int_{\partial M}f(\zeta)^{(i,k-i)}\overline{f(\zeta)}^{(i,k-i)}\,d\sigma_M(\zeta) \\ \nonumber
                                                              &=\sum_{i=0}^k {{k}\choose{i}}^2f(\omega)^{(i,k-i)}\overline{f(\omega)}^{(i,k-i)} \cdot \frac{i!(k-i)!}{N(1+k)!} \\ \nonumber
                                                              &=\frac{N}{1+k}\sum_{i=0}^k{{k}\choose{i}}f(\omega)^{(i,k-i)}\overline{f(\omega)}^{(i,k-i)} \\ \nonumber
                                                              &=\frac{N}{1+k}\sum_{i=0}^k{{k}\choose{i}}(|f_1(\omega)|^2)^i(|f_2(\omega)|^2)^{k-i} \\ \nonumber
                                                              &=\frac{N}{1+k}(|f_1(\omega)|^2+|f_2(\omega)|^2)^k \\ \nonumber
                                                              &=\frac{N}{1+k}.
\end{align}

Therefore, by the definition of $Q_t \circ f$, we have estimations
\begin{align}
\int_0^1\,dt\int_{\partial M}|Q_t \circ f(\zeta)|^2\,d\sigma_M(\zeta)&=\sum_{j=1}^K\int_{\partial M} \big|\langle f(\zeta), f(\omega_j) \rangle\big|^{2k}\,d\sigma_M(\zeta) \\ \nonumber
                                                           &= \frac{KN}{(1+k)} \\ \nonumber
                                                           &\geqslant \frac{C_4N}{4(1+k)r^2} \\ \nonumber
                                                           &>C_5.
\end{align}
This implies that there exists $t \in [0,1]$ such that
\begin{equation}
\int_{\partial M}|Q_t \circ f(\zeta)|^2\,d\sigma_M(\zeta) > C_5.
\end{equation}
Our focus will be on this $Q_t \circ f$.

\subsection{Upper Bounds}
Fix $\zeta \in \partial M$, for $m=0,1,2,\ldots$, define
\begin{displaymath}
H_m=\{\omega_j: mr \leqslant d_M(\zeta, \omega_j)<(m+1)r\}.
\end{displaymath}
If $\omega_j \in H_m$, we have
\begin{equation}
\big|\langle f(\zeta), f(\omega_j) \rangle \big|^2 \leqslant 1-m^2r^2,
\end{equation}
and $\omega_j \in E_M(\zeta, (m+1)r)$, which implies that $E_M(\omega_j, r) \subset E_M(\zeta, (m+2)r)$, 
and therefore
\begin{equation}
\sigma_M(E_M(\omega_j, r)\cdot \#H_m \leqslant \sigma_M( E_M(\zeta, (m+2)r),
\end{equation}
where $\#H_m$ denotes the cardinality of $H_m$. According to the estimation \eqref{boundaryballarea}, 
\begin{equation}
\#H_m \leqslant (m+2)^2.
\end{equation}

Thus, we have estimations
\begin{align}
|Q_t \circ f(\zeta)|&\leqslant\sum_{j=1}^K \big|\langle f(\zeta), f(\omega_j) \rangle \big|^k \\ \nonumber
&=\sum_{m=0}^{\infty}\sum_{\omega_j \in H_m}\big|\langle f(\zeta), f(\omega_j) \rangle\big|^k \\ \nonumber
&\leqslant \sum_{m=0}^{\infty} (m+2)^2(1-m^2r^2)^\frac{k}{2} \\ \nonumber
&< \sum_{m=0}^{\infty}(m+2)^2e^{-\frac{k}{2}m^2r^2} \\ \nonumber
&< \sum_{m=0}^{\infty}(m+2)^2e^{-\frac{m^2}{2}}.
\end{align}

The last series is convergent. Using $\Sigma$ to denote the sum, we therefore have
\begin{equation}
\Big|\frac{Q_t \circ f(\zeta)}{\Sigma}\Big| < 1.
\end{equation}

\subsection{RW-Sequences}

According to the results above, we define, for $\zeta \in \partial M$,
\begin{equation}
W_k \circ f(\zeta)=\frac{Q_t \circ f(\zeta)}{\Sigma},
\end{equation}
which is a polynomial of $f_1(z)$ and $f_2(z)$. 
Moreover, $W_k\big(\lambda f_1(z), \lambda f_2(z)\big)= \lambda^k W_k\circ f$.
We call this a homogeneous $f$-polynomial of degree $k$.
This leads us to the following theorem.
\begin{theorem}
There exists a positive constant $c$ such that for $k \in \mathbb{N}$ and $W_k \circ f$ as defined, we have

1. $W_k \circ f$ is an homogeneous $f$-polynomial of degree $k$,

2. $|W_k \circ f(\zeta)|\leqslant 1$, and 

3. $\int_{\partial M}|W_k \circ f|^2\,d\sigma \geqslant c$.
\end{theorem}

Letting $\mathscr{U}$ be a compact subgroup of $U(2)$, we have that
\begin{equation}
\int_{\mathscr{U}}|W_t \circ U \circ f|^2\,dU=\int_{\partial M}|W_t \circ f|^2\,d\sigma \geqslant c.
\end{equation}

If $\mu$ is a positive Borel measure on $\partial M$, we have
\begin{equation}
\int_{\mathscr{U}}dU \int_{\partial M}|W_t \circ U \circ f|^2\,d\mu \geqslant c \int_{\partial M}d\mu,
\end{equation}
and therefore we can find $U_k$ such that 
\begin{equation}
\int_{\partial M}|W_t \circ U_k \circ f|^2\,d\mu \geqslant c\int_{\partial M}d\mu.
\end{equation}

Note that the results in the theorem remain true if we consider $W_t \circ U_k \circ f$ instead of $W_t \circ f$, we  add another property for this $f$-polynomial:
\begin{proposition}
\label{EW-sequences}
There exists a positive constant $c$ such that for $k \in \mathbb{N}$ and $W_k \circ f$ as defined, we have

1. $W_k \circ f$ is an homogeneous $f$-polynomial of degree $k$,

2. $|W_k \circ f(\zeta)|\leqslant 1$, 

3. $\int_{\partial M}|W_k \circ f|^2\,d\sigma \geqslant c$, and 

4. $\int_{\partial M}|W_k \circ f|^2\,d\mu \geqslant c \int_{\partial M}d\mu$, if $\mu$ is a positive Borel measure.
\end{proposition}

\section{Inner Functions}
Any holomorphic function can be written as a series of homogeneous polynomials:
\begin{equation}
h(z)=\sum_{k=0}^{\infty}h_k(z)
\end{equation}
Then $h\circ f$ can be written as a series of $f$-homogeneous polynomials
\begin{equation}
h\circ f(z)=\sum_{k=0}^{\infty}h_k \circ f(z)
\end{equation}

Let $E$ be a set of nonnegative integers, an $(E,f)$-polynomial is a finite sum of the form $\sum_{k \in E}F_k$,
where $F_k$ is a $f$-homogeneous polynomial of degree $k$.
If  $k$ is taken from 0 to $\infty$. We call it an $(E,f)$-function.

If, additionally, $E$ is such that  there are such integers $a_m$ $(m=1, 2, 3, \ldots)$  that $E$ contains $j+a_m$ for $j=1, 2, \ldots, m$, we say that $E$ is an LI-set, which means $E$ contains arbitrarily long intervals of consecutive integers.
A quick result is that the removal of any finite subset of an LI-set still gives an LI-set.

\begin{proposition}
\label{polyapproximation}
Suppose that $\varphi \in C(\overline{B})$, $E$ is an LI-set, and for $k=1, 2, \ldots$, 
$f_k$ is an $f$-homogeneous polynomial of degree $k$, with $|f_k|\leqslant 1$ on $\overline{M}$.

Then there is a sequence $\{k_i\}$, and there are $(E,f)$-polynomials $F_i$ such that
\begin{equation}
\lim_{i \to \infty}|F_i(z)-f_{k_i}(z) \varphi\big(f(z)\big)|=0
\end{equation}
uniformly on $\overline{M}$.
\end{proposition}

The following gives a proof to this proposition, which is analogous to Rudin's in \cite{R2}, with a few necessary modifications, the most important of which is to give a Cauchy integral formula for our case.

Using the change of variables $z=f(w)$ and applying the Cauchy integral over $S$, we have
\begin{equation}
F(f(w))=\frac{1}{N}\int_{\partial M}\frac{F(f(\zeta))\,d\sigma_M(\zeta)}{(1-\langle f(w), f(\zeta) \rangle)^2}.
\end{equation}

We may assume $\varphi$ is actually a polynomial of $z_1, z_2, \overline{z}_1$ and $\overline{z}_2$, and we first consider the case $\varphi=\psi=z^{(\alpha_1,\alpha_2)}\overline{z}^{(\beta_1, \beta_2)}$.

We then define
\begin{equation}
\label{pk}
P_k=\frac{1}{N}\int_{\partial M}\frac{f_k(\zeta)\varphi\big(f(\zeta)\big)\,d\sigma_M(\zeta)}{\big(1-\langle f(z), f(\zeta) \rangle\big)^2}.
\end{equation}

Note that $f_k(\zeta)$ is a finite sum of terms of the form $f(\zeta)^{(i,k-i)}$.
We expand $(1-\langle f(z), f(\zeta) \rangle)^{-2}$, and due to \eqref{zerointegral}, we are only interested in the following integral, with the integer $t$ to be fixed:
\begin{equation}
\int_{\partial M} f(\zeta)^{(\alpha_1,\alpha_2)}\overline{f(\zeta)}^{(\beta_1, \beta_2)}f(\zeta)^{(i,k-i)}\Big(f_1(z_1) \overline{f_1(\zeta)}+f_2(z_2) \overline{f_2(\zeta_2}\Big)^t \, d\sigma_M(\zeta).
\end{equation}

We compute
\begin{equation}
\Big(f_1(z) \overline{f_1(\zeta)}+f_2(z) \overline{f_2(\zeta)}\Big)^t=\sum_{j=0}^t{{t}\choose{j}} f(z)^{(j,t-j)}\overline{f(\zeta)}^{(j,t-j)}.
\end{equation}
Thus, in order for the integral to be nonzero, we must have relations
\begin{equation}
\left\{
\begin{array}{lr}
\alpha_1+i=\beta_1+j,&\\
\alpha_2+k-i=\beta_2+t-j,& \\
\end{array} \right.
\end{equation}
from which we obtain
\begin{equation}
t=k+(\alpha_1+\alpha_2)-(\beta_1+\beta_2).
\end{equation}

So we can conclude that \eqref{pk} is actually an $f$-homogeneous polynomial of degree $k+(\alpha_1+\alpha_2)-(\beta_1+\beta_2)$.

Taking into consideration that $\varphi\circ f(z)=\sum C_{\alpha \beta} f(z)^{\alpha} \overline{f(z)}^{\beta}$, 
and that $E$ is an LI-set, we can check that $P_k$ is an $(E,f)$-polynomial for infinitely many $k$.

To check for uniform convergence, we first note that
\begin{equation}
P_k(z)-f_k(z)\varphi\big(f(z)\big)=\Xi_k(z)= \frac{1}{N}\int_{\partial M}\frac{f_k(\zeta)\Big(\varphi\big(f(\zeta)\big)-\varphi\big(f(z)\big)\Big)\,d\sigma_M(\zeta)}{\big(1-\langle f(z), f(\zeta) \rangle\big)^2}
\end{equation}
for $z \in \overline{M}, \zeta \in \partial M - \{z\}$.

Then we can check that 
\begin{equation}
\lim_{k \to \infty} \Xi_k(z) =0.
\end{equation}
Noting that $\{\Xi_k\}$ is equicontinuous completes our verification of the proposition.

Now, suppose that $\varphi$ is a positive LSC function on $\overline{B}$ and that $\varphi \circ f \in L^2(\mu)$. In fact, we may simply assume that $\varphi \in C(\overline{B})$ because we can use increasing sequences of positive continuous functions to approximate LSC $\varphi$ from below.
According to the result in Proposition \eqref{EW-sequences}, for \textit{any} $k$,
we can find $f$-homogeneous polynomials $W_k \circ f$ such that $|W_k \circ f| \leqslant 1$ and 
\begin{equation}
\int_{\partial M}|(W_k \circ f)(\varphi \circ f)|^2\,d\mu \geqslant c \int_{\partial M}(\varphi \circ f)^2\,d\mu.
\end{equation}

According to Proposition \eqref{polyapproximation}, 
we can find an $(E,f)$-polynomial $F$ and \textit{one} $W_k\circ f$ with the relation
\begin{equation}
\label{gaiainequality}
|F-(W_k \circ f)(\varphi \circ f)| < \frac{\varphi \circ f}{4}.
\end{equation}

Now, let $P=\frac{4}{5}F$. From the inequality \eqref{gaiainequality} and the fact that $|W_k \circ f| \leqslant 1$, we immediately have the result
\begin{equation}
|P|<\frac{4}{5}(1+\frac{1}{4})\varphi \circ f = \varphi \circ f.
\end{equation}

Additionally, the inequality \eqref{gaiainequality} gives us
\begin{equation}
|F|>|(W_k \circ f)(\varphi \circ f )|-\frac{\varphi \circ f}{4}.
\end{equation}
Consequently, we can see that
\begin{align}
|F|^2&=|(W_k \circ f)(\varphi \circ f )|^2-\frac{|(W_k \circ f)(\varphi \circ f )|\varphi \circ f}{2}+\frac{(\varphi \circ f)^2}{16} \\ \nonumber
&>|(W_k \circ f)(\varphi \circ f )|^2-\frac{(\varphi \circ f)^2}{2}.
\end{align}
This gives us another estimation:
\begin{align}
\int_{\partial M}|P|^2&=\frac{16}{25}\int_{\partial M}|F|^2 \\ \nonumber
                      &>\frac{16}{25}\int_{\partial M}|(W_k \circ f)(\varphi \circ f )|^2-\frac{(\varphi \circ f)^2}{2} \\ \nonumber
                      &\geqslant \frac{16}{25}(1-\frac{1}{2})\int_{\partial M}(\varphi \circ f)^2 \\ \nonumber
                      &=\frac{8}{25}\int_{\partial M}(\varphi \circ f)^2.
\end{align}

To summarize, we re-state the result in a theorem:
\begin{theorem}
\label{generatingpolynomial}
If $\varphi$ is a positive LSC function on $\overline{B}$ and $\varphi \circ f \in L^2(\mu)$, we can find an $(E,f)$-polynomial $P$ such that $|P(z)|<\varphi \big(f(z)\big)$ and $$\int_{\partial M}|P|^2\,d\mu >\frac{8}{25}\int_{\partial M}(\varphi \circ f)^2\,d\mu.$$
\end{theorem}

We can now verify the existence of inner functions.
\begin{theorem}
\label{existencethm}
\label{existence}
Suppose $\varphi >0$ on $S$, $\varphi \circ f \in LSC \cap L^2(\sigma_M)$, and $E$ is an LI-set. Then there is a nonconstant $E$-function $F \in H^2(M)$ whose boundary values $F^*$ satisfy
\begin{equation}
|F^*(\zeta)|=\varphi \circ f(\zeta) \ a.e. [\sigma_M].
\end{equation}
\end{theorem}
\begin{proof}
We start with an $(E,f)$-polynomial $P_0$ satisfying $|P_0|<\varphi \circ f$,and denote the set of integers which are degrees of monomials in $P_0$ by $E_0$.
For example, we have $P_0=0$.

If we let $Q_0=P_0$, then on $\partial M$ we have
\begin{equation}
|Q_0|<\varphi \circ f.
\end{equation}

Thus, according to Theorem \eqref{generatingpolynomial} we can construct an $(E \setminus E_0,f)$-polynomial $P_1$ such that $|P_1|< \varphi \circ f -|Q_0|$ on $\partial M$, and
\begin{equation}
\int_{\partial M}|P_1|^2\,d\sigma_M >\frac{8}{25}\int_{\partial M}(\varphi \circ f -|Q_0|)^2\,d\sigma_M.
\end{equation}
Denote the set of integers which are degrees of monomials in $P_1$ by $E_1$, and note that $E_0 \cap E_1 = \emptyset$ and therefore, $P_0$ and $P_1$ are orthogonal to each other. 
We can still regard $P_1$ as an $(E,f)$-polynomial.

Now, suppose we have found pairwise orthogonal $(E,f)$-polynomials $P_0, P_1, \ldots, P_N$ such that
\begin{equation}
|Q_N|<\varphi \circ f
\end{equation}
and
\begin{equation}
\int_{\partial M}|P_N|^2\,d\sigma_M >\frac{8}{25}\int_{\partial M}(\varphi \circ f -|Q_{N-1}|)^2\,d\sigma_M,
\end{equation}
where we define $Q_j =\sum_{i=0}^j P_i$.

Then by Theorem \eqref{generatingpolynomial} we can find an $\big(E \setminus \cup_{i=0}^{N}E_i, f\big)$-polynomial $P_{N+1}$, 
which can be also regarded as an $(E,f)$-polynomial, satisfying
\begin{equation}
\label{PandQ}
|P_{N+1}|< \varphi \circ f - |Q_N|,
\end{equation}
and
\begin{equation}
\int_{\partial M}|P_{N+1}|^2\,d\sigma_M >\frac{8}{25}\int_{\partial M}(\varphi \circ f -|Q_N|)^2\,d\sigma_M.
\end{equation}

We note that $P_{N+1}$ is orthogonal to $P_1, \ldots, P_N$. 
By definition of $Q_{N+1}$ and the inequality \eqref{PandQ}, we see that
\begin{equation}
|Q_{N+1}|\leqslant |Q_N|+|P_{N+1}|<|Q_N|+\varphi \circ f - |Q_N|= \varphi \circ f.
\end{equation}

In short, we start with an $(E,f)$-polynomial $P_0$ and then construct $Q_0$. If we already have pairwise orthogonal $(E,f)$-polynomials $P_1, \ldots, P_N$, with $Q_N=P_0+\cdots+P_N$, satisfying
\begin{equation}
\label{Qupbound}
|Q_N|<\varphi \circ f
\end{equation}
and
\begin{equation}
\label{intPlowerbound}
\int_{\partial M}|P_N|^2\,d\sigma_M >\frac{8}{25}\int_{\partial M}(\varphi \circ f -|Q_{N-1}|)^2\,d\sigma_M,
\end{equation}
then we can inductively construct an $(E,f)$-polynomial $P_{N+1}$ orthogonal to $P_i, 0 \leqslant i \leqslant N$, and $Q_{N+1}=Q_N+P_{N+1}$, such that 
\begin{equation}
|Q_{N+1}| < \varphi \circ f
\end{equation}
and
\begin{equation}
\int_{\partial M}|P_{N+1}|^2\,d\sigma_M >\frac{8}{25}\int_{\partial M}(\varphi \circ f -|Q_N|)^2\,d\sigma_M.
\end{equation}

Next, we notice that for any $N$, because of the orthogonality of $\{P_i\}$,
\begin{equation}
\label{intQupbound}
\int_{\partial M}|Q_N|^2\,d\sigma_M=\int_{\partial M}\Big|\sum_{i=0}^N P_i\Big|^2\,d\sigma_M=\sum_{i=0}^N\int_{\partial M}|P_i|^2\,d\sigma_M.
\end{equation}
On the other hand, by the inequality \eqref{Qupbound}, we have
\begin{equation}
\int_{\partial M}|Q_N|^2\,d\sigma_M <\int_{\partial M}(\varphi \circ f)^2\,d\sigma_M.
\end{equation}
Therefore we have the relation
\begin{equation}
\sum_{i=0}^{\infty}\int_{\partial M}|P_i|^2\,d\sigma_M \leqslant \int_{\partial M}(\varphi \circ f)^2\,d\sigma_M,
\end{equation}
and it makes sense to define $F = \sum_{i=0}^{\infty} P_i$, because this series converges and further more we know that $F \in H^2(M)$.

According to \eqref{intQupbound}, we can deduce that $|Q_N| \to |F^*|$ in $L^2(\sigma_M)$. However, by \eqref{intPlowerbound}, whose left-hand side goes to 0 as $N \to \infty$, we have $|Q_N| \to \varphi \circ f$ in $L^2(\sigma_M)$. Thus, we can conclude that $|F^*(\zeta)|=\varphi\circ f(\zeta) \ a.e. [\sigma_M]$, and the theorem is proved.
\end{proof}

We are finally ready to state the main result, simply a special case of Theorem \eqref{existence}, taking $\varphi \equiv 1$:

\begin{proposition}
\label{existencepro}
Inner functions exist for the domain $M$.
\end{proposition}

\section{Concluding Remark}
Our result shows the existence of inner functions a more general type of domains. 
This also provides insight for some domains of finite type. 
However, not all domains of finite type can be related to the unit ball in such a way. 
Therefore, it is still unknown whether for all domains of finite type inner functions exist. 
There are other methods and tools we may take advantage of. 
We may also take into consideration that the weakly pseudoconvex points form a set of measure zero on the boundary of domains of finite type, which is a result from the work of D.~Catlin in \cite{C}.

\bibliographystyle{amsplain}

\end{document}